\documentclass{article}
\usepackage[english]{babel}
\usepackage[T1]{fontenc}
\usepackage[latin1]{inputenc}
\usepackage[nice]{nicefrac}
\usepackage{amsmath}
\usepackage{amssymb}
\usepackage{cases}
\newtheorem{theo}{Theorem}
\newtheorem{definition}{Definition}
\newtheorem{lemme}{Lemma}
\newtheorem{rque}{Remark}

\def\dsp{\displaystyle}
\def\abs#1{\vert #1 \vert}
\def\uu{\textbf{u}}
\def\UU{\textbf{U}}
\def\kk{\textbf{k}}

\newcommand{\findem}{\begin{flushright} \vspace*{-0.4cm} $\blacksquare$ \end{flushright}}
\newcommand{\vecdeux}[2]{\left(\begin{array}{c} #1\\#2 \end{array}\right)}

\everymath{\displaystyle}
\begin{document}
\title{Compressible primitive equation: formal derivation and stability of weak solutions}

\author{Mehmet Ersoy\thanks{e-mail: Mehmet.Ersoy@univ-savoie.fr, LAMA, UMR 5127 CNRS, Universit\'e de Savoie,
73376 Le Bourget du lac  cedex, France.}, 
Timack Ngom\thanks{email: leontimack@yahoo.fr, Laboratoire d'Analyse Num\'erique et Informatique(LANI),
Universit\'e Gaston Berger de Saint-Louis, 
UFR SAT, BP 234, Saint-Louis S\'en\'egal}~and Mamadou Sy\thanks{email: syndioum@yahoo.fr, Laboratoire d'Analyse Num\'erique et Informatique(LANI),
Universit\'e Gaston Berger de Saint-Louis, 
UFR SAT, BP 234, Saint-Louis S\'en\'egal}\\[0.3cm]
Laboratoire de Math\'ematiques, Universit\'e de Savoie\\
73376 Le Bourget du Lac, France
}
\maketitle
\begin{abstract}
We present a formal derivation of a simplified version of Compressible Primitive Equations (CPEs) for atmosphere
modeling. They are obtained from  $3$-D compressible Navier-Stokes equations with an \emph{anisotropic  viscous  stress tensor} where viscosity depends on the density. We then study  the stability of the weak solutions
of this model by using an intermediate model, called model problem, which is more simple and practical, to achieve
the main result.
\end{abstract}
\textbf{Keywords}: 
Compressible primitive equations, Compressible viscous fluid,  \emph{a priori} estimates,
Stability of weak solutions.
\tableofcontents
%%%%%%%%%%%%%%%%%%%%%%%%%%%%%%%%%%%%%%%%%%%%%%%%%%%%%%%%%%%%%%%%%%%%%%
% INTRODUCTION
%%%%%%%%%%%%%%%%%%%%%%%%%%%%%%%%%%%%%%%%%%%%%%%%%%%%%%%%%%%%%%%%%%%%%%
\section{Introduction}\label{SectionIntro}
Among equations of geophysical fluid dynamics (see \cite{B81}), the equations governing the motion of the
atmosphere are the Primitive Equations (PEs). In the hierarchy of
geophysical fluid dynamics models, they are  situated between non hydrostatic models and shallow water models.
They are obtained from  the full $3$ dimensional set of Navier-Stokes equations  for atmosphere modeling,
\begin{eqnarray}
\rho \frac{D}{Dt}\UU +  \nabla p + \rho \textbf{g}=
D,\label{FullSetAtmosphereModeling1}\\
\frac{D}{Dt}\rho +\rho \mbox{div} \UU  =0,\label{FullSetAtmosphereModeling2}\\
c_p \frac{D}{Dt}T -\frac{1}{T}\frac{D}{Dt}p = Q_T,\label{FullSetAtmosphereModeling3}\\
\frac{D}{Dt} q = Q_q,\label{FullSetAtmosphereModeling4}\\
p =  RT\,\rho \label{FullSetAtmosphereModeling5}
\end{eqnarray} where 
$$
\frac{D}{Dt} = \partial_t + \UU\cdot \nabla \,.
$$
$\UU$ is the three dimensional velocity vector with component $\uu$ for horizontal velocity and $v$ for the
vertical one. The terms  $\rho$, $p$, $T$, $\textbf{g}$ stand for the density, the pressure, the temperature, the gravity
vector $(0,0,g)$. The diffusion term $D$ is written as:
\begin{equation}\label{TemamViscosity}
D= \mu \Delta_x \UU + \nu \partial_{yy}^2 \UU
\end{equation}
where $ \Delta_x $ stands for the derivatives of second order with respect to  the horizontal variables
$x=(x_1,x_2)$, and 
$\mu \neq \nu$ represents the anisotrope pair of viscosity. 
The  diffusive term $Q_q$  represent the molecular diffusion  where $q$ is the amount of water in the air,  and $Q_T$ is the heat diffusion standing for the solar heating (see for
instance \cite{TZ04} for details of diffusive terms). The last
term $c_p$ is the specific heat of the air at constant pressure and $R$ is the specific gas constant for the air.

A scale analysis show that only the terms $\partial_y p$ and $g
\rho$ are dominant (see e.g. \cite{Pedlowski87}). This leads to replace the third equation of 
(\ref{FullSetAtmosphereModeling1}) with the hydrostatic one to obtain the so-called Compressible Primitive
Equations (CPEs) for atmosphere modeling, 
\begin{equation}\label{CPEGeneral}
\left\{
\begin{array}{l}
\rho \frac{d}{dt}\uu +  \nabla_x p = D,\\
\partial_y p = - g \rho, \\
\frac{d}{dt}\rho + \rho \mbox{div} \UU =0,\\
c_p \frac{D}{Dt}T -\frac{1}{T}\frac{D}{Dt}p = Q_T,\\
\frac{D}{Dt} q = Q_q,\\
p =  RT\,\rho 
\end{array}
\right.
\end{equation} where $x$ and $y$ stand for the horizontal and vertical coordinate and
$$
\frac{d}{dt} = \partial_t + \uu\cdot \nabla_x + v \partial_y\,.
$$
\paragraph{\textbf{Derivation ``background''}}
In this paper, we present  the derivation of  Compressible  Primitive 
Equations (CPEs) close to Equations (\ref{CPEGeneral}) (without taking in account complex phenomena such as 
the amount of water in the air and  the solar heating) from the
$3$-D Navier-Stokes equations 
with an \emph{anisotropic viscous tensor}.  Emphasizing to the difference of sizes of the vertical and horizontal
dimensions in the atmosphere (10 to 20 for  height with respect to thousands of kilometers of length),   
we derive the hydrostatic balance approximation for the vertical motion. We obtain  \emph{simplified CPEs}:
\begin{equation}\label{Model2}
\left\{\begin{array}{l}
\dsp \partial_{ t} {\rho} + \mbox{div}_{x}\left({\rho}\,{\uu} \right) 
+ \partial_{y}\left(\rho {v}\right) =  0,\\ 
\dsp \partial_{ t} \left({\rho}\,{\uu}\right) +  \mbox{div}_{{x}}\left({\rho}\,\uu \otimes
{\uu}\right) + \partial_{{y}}\left({\rho}\,{v}{\uu}\right) +
\nabla_{x} p(\rho) = \mbox{div}_{x}\left(\nu_1 D_{{x}}({\uu})\right)\\ 
\dsp  + \partial_{y} \left( \nu_2 \partial_{{y}}{\uu}\right), \\ 
\dsp \partial_y p(\rho) = -g \rho 
\end{array}\right.
\end{equation} where $y$ stands for  the vertical coordinate. The main difference between Model (\ref{CPEGeneral})
and Model (\ref{Model2}) is the viscous term.
Moreover, if $p= c^2  \rho$ with $c=RT$ (for instance, as above), the density $\rho$ is written $ \xi(t,x)e^{-g/c^2
y}$ where $\xi$ , called again ``density'', is an unknown of the following system called \emph{model problem}:
\begin{equation}\label{Model3}
\left\{\begin{array}{l}
\dsp \frac{d}{dt}\xi + \xi \mbox{div}_x(\uu) + \xi \partial_z w=0,\\ 
\dsp \xi \frac{d}{dt}\uu  + c^2\nabla_{x} (\xi) =  D, \\
\dsp \partial_z \xi = 0  
\end{array}\right.
\end{equation}
which is obtained from System (\ref{Model2}) by the ``simple'' change of variables  
$z = 1-e^{-y}$ and $w=e^{-y} v$ where  $$
\frac{d}{dt} = \partial_t + \uu\cdot \nabla_x + w \partial_z
$$ and $D$ stands for the following viscous terms
\begin{equation}\label{Viscousness}
D= \mbox{div}_{x}\left(\nu_1 D_{{x}}({u})\right) + \partial_{z} \left( \nu_2
\partial_{{z}}{u}\right).
\end{equation}
As explained below, we cannot obtain a result directly on System (\ref{Model2}) without using the intermediate Model (\ref{Model3}), so we use the fact that System (\ref{Model3}) is very close to System (\ref{Model2}), and that the equation $\partial_z \xi=0$ is one of the
key ingredient to achieve the stability of weak solutions of System (\ref{Model3}), to propagate the result to System (\ref{Model2}).
% 
% The model problem (\ref{Model3}), already introduced in \cite{K36} and \cite{GK05,EN10_3} for the $2$ dimensional
% version, from its viscous form it is also close to the shallow water equations, energetically consistent
% \cite{G93}, since it may be easily derived from this one by vertical averaging and the
% use of appropriate barotropic law \cite{EN10_2}. 
% Thus Equations (\ref{Model3}) can be seen as an intermediate
% model from the general CPEs and the shallow water equations.
\paragraph{\textbf{Mathematical ``background''}}
The mathematical study of PEs for atmosphere modeling were first studied by J.L. Lions, R. Temam and S. Wang (\cite{LTW92}) where they produced the mathematical formulation  of System (\ref{CPEGeneral}) in $2$ and $3$
dimensions based on the works 
of J. Leray and they obtained the existence of weak solutions for all time (see also \cite{TZ04} where the  result was proved by different means). For instance, in \cite{TZ04},  using the
hydrostatic equation, they used the pressure $p$ as vertical coordinate instead of the altitude $y$. Moreover, they
wrote System (\ref{CPEGeneral}) in spherical coordinates $(\phi,\theta,p)$ to change compressible equations
to incompressible ones to  use the well-known results of incompressible theory. They distinguished the
\emph{prognostic} variables from the \emph{diagnostic} variables, which are: $(\uu,T,q)$
for the prognostic and $(v,\rho,\Phi)$ for the diagnostic variables where $\Phi$ is the geopotential $g
y(\phi,\theta,p,t)$. Diagnostic variables $(v=v(\uu),\rho=\rho(T),\Phi=\Phi(T))$ can be written as a
function of the prognostic variables through the div-free equation, $p=RT \rho$ and by integrating the mass
equation which is written in the new coordinates as follows:$$\partial_p \Phi + \frac{RT}{p}=0.$$
Then, the outline of their proof of the existence was: they wrote
\begin{itemize}
\item a weak formulation of the PEs (by defining appropriate space functions) of the form $$\frac{d U}{dt} +
AU+B(U,U)+E(U)=l$$ where $U=(\uu,T,q)$ with initial data $U(0)=U_0$ and $A,B,E$ are appropriate functional,
\item finite differences in time: $U^n$,
\item \emph{a priori estimates} for $U^n$,
\item approximate functions: $U_{\Delta t}(t) = U^n $ on $((n-1)\Delta t,n \Delta t)$ (following \cite{Temam77})
\item \emph{a priori estimates} for $U_{\Delta t}$,
\end{itemize}
and they proved the passage to the limit.

Setting $T$ and $q$ constant, the main difference between Model (\ref{CPEGeneral}) and (\ref{Model2}) comes from
the viscous term (\ref{TemamViscosity}) and (\ref{Viscousness}). Starting from the Navier-Stokes equations with
\emph{non constant density dependent viscosity} and \emph{anisotropic viscous tensor}, it is natural to get the
viscous term (\ref{Viscousness}). This term is also present in the viscous shallow water equations (see e.g.
\cite{HandBookBresch09}).  
% The global existence and uniqueness of solution of shallow water equations or CPEs equations 
% is one of the most challenging mathematical problems in applied analysis.
% Consequently, the simplified model (\ref{Model2}) which is an intermediate model 
% is equally challenging from this point
% of view.
In the same spirit   than \cite{TZ04}, authors \cite{EN10_3} showed  a global weak
existence for the $2$-D version of model problem (\ref{Model2})  with $p(\rho) =c^2 \rho$ (with the notation
above) by a change of vertical coordinates (but not on $p$ as done in \cite{TZ04}) which led to  prove that Model
(\ref{Model2}) and Model
(\ref{Model3}) are equivalent. Then, using  an existence result provided by Gatapov  \emph{et al}
\cite{GK05} for Model (\ref{Model3}), they  could conclude.  Existence result \cite{GK05} for System (\ref{Model3})
was obtained as follows:
\begin{itemize}
 \item by a useful change of variables in the
Lagrangian coordinates, authors \cite{GK05} showed that the density $\xi$ is bounded from above and below . 
\item  by \emph{a priori} estimates and writing the system for the oscillatory part of the velocity $\uu$, authors
\cite{GK05} obtained the existence result thanks to a Schauder fix point theorem.
\end{itemize}
Unfortunately, this approach \cite{GK05} fails for the $3$-D version (\ref{Model3}) since the change of variables
in Lagrangian coordinates does not provide enough information to bound the density $\xi$. 
Moreover,  to show a stability result for weak solutions for Model (\ref{Model2}) with standard techniques
also fails, since multiplying the conservation of the momentum equations of Model (\ref{Model2}) by
$(u,v)$ gives:
$$
\frac{d}{dt}\int_{\Omega} \rho\abs{u}^2+ \rho\ln\rho-\rho + 1 \,dx dt +  
\int_{\Omega} \nu_1(\rho) \abs{D_x(u)}^2 + nu_2(\rho)\abs{\partial_{yy}^2 u} \,dx  +  
\int_{\Omega} \rho g v \,dx
$$
where the sign of the integral $\dsp \int_{\Omega} \rho g v \,dx$ is unknown  (in the equation above, $D_x(u)$
stands for $\dsp \frac{\nabla_x u +\nabla_x^t  u }{2}$). It appears that, \emph{prima facie}, there is a missing information  on $v$ to avoid the integral term $\dsp \int_{\Omega} \rho g v \,dx$ introduced by the hydrostatic
equation $\partial_y p =-g \rho$. In fact, the study of the weak  solutions stability   cannot be performed directly on System (\ref{Model2}) (at least up to our knownledge): therefore,  we have to study  the intermediate Model (\ref{Model3}) through the change of vertical coordinates.
Indeed, we have just remarked that (as described above) that $\rho$ is written as $\xi e^{-g/c^2 y}$, where $\xi$ does
not depend on $y$. Thus, performing a change
of variable in vertical coordinate in Model (\ref{Model2}), following \cite{EN10_3}, we showed that Model
(\ref{Model2}) could be written as
Model (\ref{Model3}). 
As the hydrostatic equation in System (\ref{Model3}) is $\partial_z \xi = 0$,
an energy equality is easily
obtained, and provided some \emph{a priori} estimates. Nevertheless, those estimates are not strong enough  to
pass to the limit in the non linear terms; additional informations are required. On the other hand, we have to
remark that the missing information for the
vertical speed $v$ for Model (\ref{Model2}) (or equivalently $w$ for Model
(\ref{Model3})) is fulfilled by  the equation of the mass of System (\ref{Model3}), which is also written as:
$$\partial_{zz}^2 w = \frac{1}{\xi}\mbox{div}_x(\xi\partial_y \uu)\,.$$ 
Then, the fact that $\xi=\xi(t,x)$ combined with the equation above, allowed to obtain a mathematical entropy, the
BD-entropy (initially introduced in \cite{BD04},
where a simple proof was given in \cite{BDGG07,BDGV07} or in \cite{HandBookBresch09} and the reference therein). 
Let us also notice that, as for shallow water equations (see e.g. \cite{BD02,BD03,BD06} to cite only a few),  it is
necessary to add a regularizing  term (as capillarity of friction) to equations (\ref{Model3}) (equivalently to Model (\ref{Model2})) to conclude to the stability of weak solutions for Model (\ref{Model3}): in this present case, we add a
quadratic friction source term which is written $r\rho u\abs{u}$ for System (\ref{Model2}) or equivalently $r\xi
u\abs{u}$ for System (\ref{Model3}). Indeed, the viscous term (\ref{Viscousness}) combined to  the friction term  brings some 
regularity on the density, which is required to pass to the limit in the non linear terms
(e.g. for the term $\rho u \otimes u$, where typically a strong convergence of  $\sqrt{\rho} u$ is needed).
Finally, by the reverse change of variables, the estimates, necessary to prove stability of weak solutions, were
obtained for System (\ref{Model2}) from those of System (\ref{Model3}). 

We  note that, for the sake of simplicity, periodic conditions on the spatial horizontal domain
$\Omega_x$ are assumed, since it avoids an incoming boundary term (whose sign is unknown: see e.g. \cite{BDGV07}), which appears when we
seek a mathematical BD-entropy. Let us also precise that ``good'' boundary conditions on $\Omega_x$ may be used
(see \cite{BDGV07})  instead of periodic ones to avoid this boundary term. 

We may also perform this analysis without the quadratic friction term by using the ``new'' multiplier introduced in
\cite{MV07} which provides another mathematical entropy: particularly to estimate bounds of $\rho
u^2$ in a better space than $L^{\infty}(0,T;L^{1}(\Omega))$.

\paragraph{\textbf{This paper is organized as follows}} In Section \ref{SectionDerivation}, starting from the $3$-D
compressible Navier-Stokes equations 
with  an \emph{anisotropic viscous tensor}, we formally derive the simplified Model (\ref{Model2}) as described
above. We present the main result in Section \ref{SectionMainResult}. In the third and last Section
\ref{SectionMathematicalStudy}, we prove the
main result. Firstly, we  show that Model (\ref{Model2}) can be
rewritten as Model (\ref{Model3}) which is more simpler. Then, taking advantages of the property of
the density $\xi$, adding a quadratic  friction term (following \cite{BD02,BD03}), we obtain a mathematical energy
and entropy which provides enough estimates to pass to the limit in 
Model (\ref{Model3}). Finally, following \cite{EN10_3}, the stability result for Model
(\ref{Model2}) is easily obtained.
%%%%%%%%%%%%%%%%%%%%%%%%%%%%%%%%%%%%%%%%%%%%%%%%%%%%%%%%%%%%%%%%%%%%%%
% Formal derivation of the simplified atmosphere model
%%%%%%%%%%%%%%%%%%%%%%%%%%%%%%%%%%%%%%%%%%%%%%%%%%%%%%%%%%%%%%%%%%%%%%
\section{Formal derivation of the  simplified atmosphere model}\label{SectionDerivation}
We consider the Navier-Stokes model in a bounded three dimensional domain 
with periodic boundary conditions on $\Omega_x$ and free conditions on the rest of the boundary. 
More exactly,  we assume that motion of the medium occurs in a domain $\Omega =
\{(x,y);\,x\in \Omega_x,\, 0<y<h\}$
where $\Omega_x= \mathbb{T}^2$ is a torus.
The full Navier-Stokes equation is written:
\begin{eqnarray}
\partial_t\rho+\mbox{div}(\rho u)=0,\label{MassEquationCNSA}\\
\partial_t(\rho u) +\mbox{div}(\rho u\otimes u)-\mbox{div} \sigma(u) -\rho f=0, \label{MomentumEquationCNSA}\\
p=p(\rho)\label{pressure}
\end{eqnarray} 
where  $\rho$ is the density of the fluid and $u=(\uu,v)^t$ stands for the fluid velocity with  
$\uu=(u_1,u_2)^t $ the horizontal component and $v$ the vertical one. 
The pressure law is given by the  equation of state:
\begin{equation}\label{PressureLaw}
p(\rho)= c^2\rho
\end{equation} for  some given constant $c$. 
The term $f$ is the quadratic friction source term and  the gravity strength is given as follows: 
$$f=-r\sqrt{u_1^2+u_2^2}\,(u_1,u_2,0)^t-g\kk$$ where $r$ is  a positive constant
coefficient, $g$ 
is the gravitational constant and $\kk = (0,0,1)^t$ (where $X^t$ stands for the transpose of tensor $X$).
The term $\sigma(u)$ is a non symmetric stress with the following viscous tensor (see e.g. \cite{K36,GK05,EN10_3})
$\Sigma(\rho)$:
$$
\left(
\begin{array}{ccc}
\mu_1(\rho) & \mu_1(\rho) & \mu_2(\rho) \\
\mu_1(\rho) & \mu_1(\rho) & \mu_2(\rho) \\
\mu_3(\rho) & \mu_3(\rho) & \mu_3(\rho)
\end{array}\right)\,.
$$ 
The total stress tensor is written: 
$$\sigma(\uu)=-p I_3 + 2\Sigma(\rho):D(u)+ \lambda(\rho)\mbox{div}(u)\,I_3$$ 
where the term $\Sigma(\rho):D(u)$ is written:
\begin{equation}\label{AnisotropicStressTensor}
\left(
\begin{array}{cc}
2\mu_1(\rho)D_x(\uu) & \mu_2(\rho)\left( \partial_y \uu +\nabla_x v \right) \\ 
\mu_3(\rho)\left( \partial_y \uu +\nabla_x v \right)^t & 2\mu_3(\rho)\partial_y v 
\end{array}\right)
\end{equation}
with $I_3$ the identity matrix. 
The term $D_x(\uu)$ stands for the  strain tensor, that is:
$D_x(\uu) = \dsp\frac{\nabla_x \uu+\nabla_x^t \uu}{2}$ where 
$\nabla_x = \vecdeux{\partial_{x_1}}{\partial_{x_2}}$. 
\begin{rque}
Let us remark that,  if we play with  the magnitude of  viscosity $\mu_i$, the matrix $\Sigma(\rho)$ will be useful
to set a privileged flow direction.
\end{rque}
\noindent The last term $\lambda(\rho)\mbox{div}(u)$ is the classical  
normal stress tensor  with $\lambda(\rho)$ the viscosity. 
\noindent The Navier-Stokes system is closed 
with the following boundary conditions on $\partial \Omega$:
\begin{equation}\label{BoundaryConditions}
\begin{array}{l}
\textrm{ periodic conditions on }  \partial\Omega_x,\\
v_{|y=0} = v_{|y=h} = 0,\\
{\partial_y \uu}_{|y=0} = {\partial_y \uu}_{|y=h}=0.
\end{array}
\end{equation} We also assume  that the distribution of the horizontal 
component of the velocity $\uu$ and the density distribution are known at the initial time $t=0$:
\begin{equation}\label{InitialCondition}
\begin{array}{l}
 \uu(0,x,y) =\uu_0(x,y),\\
 \rho(0,x,y) =\xi_0(x) e^{-g/c^2 y}.
 \end{array}
\end{equation}
The fact that the initial condition for the density $\rho$ has the form (\ref{InitialCondition}) 
is justified at the end of Section \ref{SectionFormalDerivation}.

\noindent We assume that $\xi_0$ is a bounded positive function:
\begin{equation}\label{BoundOnXi0}
0\leqslant\xi_0(x)\leqslant M<+\infty .
\end{equation}
%%%%%%%%%%%%%%%%%%%%%%%%%%%%%%%%%%%%%%%%%%%%%%%%%%%%%%%%%%%%%%%%%%%%%%
% Formal derivation of an  atmosphere model
%%%%%%%%%%%%%%%%%%%%%%%%%%%%%%%%%%%%%%%%%%%%%%%%%%%%%%%%%%%%%%%%%%%%%%
\subsection{Formal derivation of the simplified CPEs}\label{SectionFormalDerivation}
Taking advantages of the shallowness of the atmosphere, 
we assume that the characteristic scale for the altitude  $H$  
is small with respect to the characteristic length $L$. 
So, the ratio of the vertical scale to the horizontal one is
assumed small. 
In this context, we  assume that the vertical movements 
and variations are very small compared to the horizontal ones, which justifies the following approximation.\\
Let $\varepsilon$ be a ``small'' parameter  such as:
$$\varepsilon = \frac{H}{L} = \frac{V}{U}$$ where $V$ and  $U$ are respectively the characteristic scale
of the vertical and horizontal velocity. We introduce  the characteristic time $T$ such as: $\dsp T = \frac{L}{U}$
and the pressure unit  $P = \overline{\rho}\, U^2$ where $\overline{\rho}$ is a
characteristic density. Finally, we note the dimensionless quantities of time, space, fluid
velocity, pressure, density and viscosities: 
$$\widetilde{t}=\frac{t}{T},\quad \widetilde{x}=\frac{x}{L},\quad\widetilde{y}=\frac{y}{H},\quad
\widetilde{u}=\frac{\uu}{U},\quad \widetilde{v}=\frac{v}{V},$$
$$\quad \widetilde{p}=\frac{p}{\bar{\rho}U^2},\quad
\widetilde{\rho}=\frac{\rho}{\bar{\rho}},
\quad\widetilde{\lambda}=\frac{\lambda}{\bar{\lambda}},
\quad \quad\widetilde{\mu}_j=\frac{\mu_j}{\bar{\mu}_j},
j=1,2,3$$ With these notations, the Froude number  $F_r$, the Reynolds number
associated to
the viscosity $\mu_i$ (i=1,2,3), $Re_{i}$, the Reynolds number
associated to the viscosity $\lambda$, $Re_{\lambda}$, and the Mach number  $M_a$ are written respectively:
\begin{equation}\label{FroudeReynoldsMach}
F_r= \dsp  \frac{U}{\sqrt{g\,H}},\quad 
Re_{i}  =\dsp \frac{\overline{\rho} U L}{\overline{\mu_i}}, \quad
Re_{\lambda}  =\dsp \frac{\overline{\rho} U L}{\overline{\lambda}}, \quad
M_a =\dsp  \frac{U}{c}.
\end{equation}
Applying this scaling, System (\ref{MassEquationCNSA})--(\ref{PressureLaw}) is written:
\begin{equation}\label{AdimThinLayerTheSystem}
\left\{\begin{array}{l}
\dsp \frac{1}{T}\partial_{\widetilde t} \widetilde{\rho} + \frac{U}{L} \mbox{div}_{\widetilde
x}\left(\widetilde{\rho}\,\widetilde{u} \right) + \frac{V}{H} \partial_{\widetilde{y}}
\left(\widetilde{\rho} \widetilde{v}\right) = 0,\\ 
\dsp \frac{\overline{\rho}\,U}{T}\partial_{\widetilde t} \left(\widetilde{\rho}\,\widetilde{u}\right) 
+  \frac{\overline{\rho}\,U^2}{L}\mbox{div}_{\widetilde{x}}\left(\widetilde{\rho}\,\widetilde{u}\otimes
\widetilde{u}\right)+  \frac{\overline{\rho}\,U\,V}{H} \partial_{\widetilde{y}}
\left(\widetilde{\rho}\,\widetilde{v}\widetilde{u}\right) + \frac{c^2\,\overline{\rho}}{L}
\nabla_{\widetilde{x}} \widetilde{\rho} = \\ 
\dsp \frac{\overline{\mu_1}\,U}{L^2}\mbox{div}_{\widetilde{x}}\left(\mu_1
D_{\widetilde{x}}(\widetilde{u})\right) +
\frac{\overline{\mu_2}\,U}{H^2} \partial_{\widetilde{y}} \left(\widetilde{\mu_2}\,\partial_{\widetilde{y}}
\widetilde{u}\right) +
\frac{\overline{\mu_2}\,V}{L\,H} \partial_{\widetilde{y}} \left(\widetilde{\mu_2}\,\nabla_{\widetilde{x}}
\widetilde{v}\right) + \\ 
\dsp \frac{\overline{\lambda}\,U}{L^2 } \nabla_{\widetilde{x}}
\left(\widetilde{\lambda}\,\mbox{div}_{\widetilde{x}}\left(\widetilde{u}\right) \right) + 
\dsp \frac{\overline{\lambda}\,V}{L\,H} \nabla_{\widetilde{x}}
\left(\widetilde{\lambda}\,\partial_{\widetilde{y}} \widetilde{v}\right),\\
\dsp \frac{\overline{\rho}\,V}{T}\partial_{\widetilde t} \left(\widetilde{\rho}\,\widetilde{v}\right) 
+  \frac{\overline{\rho}\,U\,V}{L}\mbox{div}_{\widetilde{x}}\left(\widetilde{\rho}\,\widetilde{u}\,
\widetilde{v}\right)+  \frac{\overline{\rho}\,V^2}{H} \partial_{\widetilde{y}}
\left(\widetilde{\rho}\,\widetilde{v}^2\right) + \frac{c^2\,\overline{\rho}}{H}
\partial_{\widetilde{y}} \widetilde{\rho} = \\
\dsp - g\,\overline{\rho} \widetilde{\rho}+ \frac{\overline{\mu_3}\,U}{LH}
\mbox{div}_{\widetilde{x}}\left(\widetilde{\mu_3}\,\partial_{\widetilde{y}} \widetilde{u}\right)        
+ \frac{\overline{\mu_3}\,V}{L^2}\mbox{div}_{\widetilde{x}}\left(\widetilde{\mu_3}\,\nabla_{\widetilde{x}}
\widetilde{v}\right)
+2\frac{\overline{\mu_3}\,V}{H^2}\partial_{\widetilde{y}}(\widetilde{\mu_3}
\partial_{\widetilde{y}}\widetilde{v
} ) \\
\dsp +\frac{\overline{\lambda}\,U}{L\,H } \partial_{\widetilde{y}}
\left(\widetilde{\lambda}\,\mbox{div}_{\widetilde{x}}\left(\widetilde{u}\right) \right) + 
\dsp \frac{\overline{\lambda}\,V}{H^2} \partial_{\widetilde{y}}
\left(\widetilde{\lambda}\,\partial_{\widetilde{y}} \widetilde{v}\right).
\end{array}\right. 
\end{equation}
Using the definition of the dimensionless number (\ref{FroudeReynoldsMach}), dropping $\widetilde{.}$,  multiplying the
mass equation of System (\ref{AdimThinLayerTheSystem}) by $T$, the momentum equation for $\uu$ 
of System
(\ref{AdimThinLayerTheSystem}) by  $\dsp
\frac{T}{\overline{\rho}\,U}$,  the momentum equation for $v$ of System
(\ref{AdimThinLayerTheSystem}) by  $\dsp
\frac{T}{\overline{\rho}\,V}$, we get the
non-dimensional version of System (\ref{MassEquationCNSA})--(\ref{PressureLaw}) as follows:
\begin{equation}\label{AdimThinLayerTheSystemWithEpsilon}
\left\{\begin{array}{l}
\dsp \partial_{ t} {\rho} + \mbox{div}_{x}\left({\rho}\,\uu\right) + \partial_{{y}}\left({\rho}
{v}\right) = 0,\\
\dsp \partial_{ t} \left({\rho}\,\uu\right) +  \mbox{div}_{{x}}\left({\rho}\,\uu\otimes
\uu\right)+  \partial_{{y}}\left({\rho}\,{v}\uu\right) +
\frac{1}{M_a^2}\nabla_{{x}} {\rho} = \frac{1}{Re_1}\mbox{div}_{{x}}\left(\mu_1 D_{{x}}(\uu)\right)
\\
\dsp  +
\frac{1}{Re_2} \partial_{{y}} \left( {\mu_2} \left(\frac{1}{\varepsilon^2}\partial_{{y}}\uu +
\nabla_{{x}}{v}\right)\right) + \frac{1}{Re_{\lambda}}
\nabla_{{x}} \left({\lambda}\,\mbox{div}_x(\uu) + \lambda\partial_y v \right),\\
\dsp \partial_{ t} \left({\rho}\,{v}\right) 
+  \mbox{div}_{{x}}\left({\rho}\,\uu\,
{v}\right)+  \partial_{{y}}
\left({\rho}\,{v}^2\right) + \frac{1}{\varepsilon^2}\frac{1}{M_a^2}\partial_{{y}} {\rho} = -
\frac{1}{\varepsilon^2}\frac{1}{F_r^2}\rho \\
\dsp+ \frac{1}{Re_3} \mbox{div}_x \left( {\mu_3} \left(\frac{1}{\varepsilon^2}\partial_{{y}}\uu +
\nabla_{{x}}{v}\right)\right)+ \frac{2}{\varepsilon^2 Re_3}\partial_{y}(\mu_3\partial_{y} v)
\\
\dsp+\frac{1}{\varepsilon^2\,Re_{\lambda}}
\partial_{{y}} \left({\lambda}\,\mbox{div}_x(\uu) + \lambda\partial_y v\right).
\end{array}\right. 
\end{equation}
Next, if we assume the following asymptotic regime:
\begin{equation}\label{RegimeAsymptotics}
\frac{\mu_1(\rho)}{Re_1} =  \nu_1(\rho),\,\frac{\mu_i(\rho)}{Re_i} =  \varepsilon^2 \nu_i(\rho),\,i=2,3 \textrm{ and }
\frac{\lambda(\rho)}{Re_{\lambda}}  = 
\varepsilon^2 \gamma(\rho).
\end{equation} and drop all terms of order $O(\varepsilon)$, 
System (\ref{AdimThinLayerTheSystemWithEpsilon}) 
reduces to the following model:
\begin{equation}\label{TheSystemWithHydrostaticApproximation}
\left\{\begin{array}{l}
\dsp \partial_{ t} {\rho} + \mbox{div}_{x}\left({\rho}\,{\uu} \right) 
+ \partial_{y}\left(\rho {v}\right) =  0,\\
\dsp \partial_{ t} \left({\rho}\,\uu\right) +  \mbox{div}_{{x}}\left({\rho}\,\uu\otimes
\uu\right) + \partial_{{y}}\left({\rho}\,{v}\uu\right) +
\frac{1}{M_a^2}\nabla_{{x}} p(\rho) = \mbox{div}_{{x}}\left(\nu_1 D_{{x}}(\uu)\right)\\ 
\dsp  +\partial_{y} \left( \nu_2 \partial_{{y}}\uu\right)+\rho f, \\
\dsp \partial_y p(\rho) = -\frac{M_a^2}{F_r^2} \rho,
\end{array}\right.
\end{equation}
called \emph{simplified CPEs.} In the sequel, we simplify by setting $M_a=F_r$. Then, the hydrostatic equation of System
(\ref{TheSystemWithHydrostaticApproximation}) with the pressure law (\ref{PressureLaw}) provides the density as 
\begin{equation}\label{DensityStratifie}
\rho(t,x,y) = \xi(t,x) e^{-y} 
\end{equation} 
for some function $\xi=\xi(t,x)$ also called ``density''. Let us note that the density $\rho$ is stratified: it
means that for any altitude $y$, the density $\rho$ has the profile of the 
function $\xi$. Therefore, Equation (\ref{DensityStratifie}) justifies 
the choice of the initial data (\ref{InitialCondition}) for the density $\rho$ at the time $t=0$.
In the sequel,  we also assume that:
\begin{equation}\label{Viscosity}
\nu_i(\rho) = \nu \rho,\,i=1,2, \textrm{ for } \nu>0.
\end{equation}
%%%%%%%%%%%%%%%%%%%%%%%%%%%%%%%%%%%%%%%%%%%%%%%%%%%%%%%%%%%%%%%%%%%%%%
% Formal derivation of an  atmosphere model
%%%%%%%%%%%%%%%%%%%%%%%%%%%%%%%%%%%%%%%%%%%%%%%%%%%%%%%%%%%%%%%%%%%%%%
\subsection{The main result}\label{SectionMainResult}
Assuming the viscosity under the form (\ref{Viscosity}) and $M_a=F_r$, we define:
\begin{definition}
A weak solution of System (\ref{TheSystemWithHydrostaticApproximation}) on $[0,T]\times\Omega$, with boundary 
(\ref{BoundaryConditions}) and 
initial conditions (\ref{InitialCondition}), is a collection  of functions $(\rho,\uu,w)$, if
$$
\begin{array}{ll}
\dsp \rho\in L^{\infty}(0,T;L^{3}({\Omega})), & \dsp \sqrt{\rho}\in L^{\infty}(0,T;H^{1}({\Omega})),\\
\dsp \sqrt{\rho}\uu\in L^2(0,T;L^2({\Omega})^2),& \dsp \sqrt{\rho} v\in L^\infty(0,T;(L^2({\Omega}))),\\
\dsp \sqrt{\rho}D_x(\uu)\in L^2(0,T;(L^2({\Omega}))^{2\times2}), &\dsp \sqrt{\rho}\partial_y v\in
L^2(0,T;L^2({\Omega})),\\
\nabla\sqrt{\rho}\in L^2(0,T;(L^2(\Omega))^3)& 
\end{array}
$$
with  $\rho\geqslant 0$ and where $(\rho,\sqrt{\rho}\uu,\sqrt{\rho} v)$ satisfies:
\begin{equation}\label{MassEquationDistribution1}
\left\{
\begin{array}{l}
\dsp\partial_t \rho + \mbox{div}_x(\sqrt{\rho} \sqrt{\rho}\uu )+ \partial_y(\sqrt{\rho}\uu \sqrt{\rho} v)=0,\\
\rho(0,x) = \rho_0(x)
\end{array}
\right.
\end{equation} in the distribution sense, 
and the following equality holds for all smooth test function $\varphi$ 
with compact support such as $\varphi(T,x,y)=0$ 
and $\varphi_0=\varphi_{t=0}$:
\begin{equation}
\begin{array}{l}
 -\int_0^T\int_{\Omega}\rho \uu\partial_t\varphi\,dx dy dt+
\int_0^T\int_{\Omega}\left(2\nu\rho D_x(\uu)-\rho \uu\otimes \uu\right):\nabla_x\varphi \,dx dy dt \\
-\int_0^T\int_{\Omega}\rho v\uu\partial_y\varphi \,dx dy dt
 -\nu\int_0^T\int_{\Omega}\rho \uu\partial_{yy}^2\varphi \,dx dy dt
+\int_0^T\int_{\Omega} r\rho\vert \uu\vert \uu\varphi\,dx dy dt \\
 -\int_0^T\int_{\Omega} \rho\mbox{div}(\varphi)\,dx dz dt + \int_0^T\int_{\Omega}  \rho v \varphi\,dx
dz dt  =\int_{\Omega} \rho_{0}\uu_{0}\varphi_0 \,dx dy.
\end{array}
\end{equation}
\end{definition} 
Now, we state the main result of this paper:
\begin{theo}\label{MainRes}
Let $(\rho_n,\uu_n,v_n)$ be a sequence of weak solutions of System (\ref{TheSystemWithHydrostaticApproximation}), with boundary 
(\ref{BoundaryConditions}) and  initial  conditions (\ref{InitialCondition}),   satisfying entropy  inequalities 
(\ref{EnergyEstimates}) and (\ref{BDEntropyEquality}) such as 
\begin{equation}
\rho_n\geqslant 0,\quad \rho_0^n \to \rho_0 \textrm{ in } L^1({\Omega})
,\quad \rho_0^n \uu_0^n \to \rho_0 \uu_0 \textrm{ in } L^1({\Omega}). 
\end{equation}
Then, up to a subsequence,
\begin{itemize}
\item $\rho_n$ converges strongly in $\mathcal{C}^0(0,T; L^{3/2}({\Omega}))$,
\item $\sqrt{\rho_n} \uu_n$ converges strongly in $L^2(0,T; L^{3/2}({\Omega})^2)$,
\item $\rho_n u_n$ converges strongly in $L^1(0,T; L^{1}({\Omega})^2)$ for all $T>0$,
\item $(\rho_n,\sqrt{\rho_n} \uu_n,\sqrt{\rho_n} v_n)$ converges to a weak solution 
of System (\ref{MassEquationDistribution1}),
\item $(\rho_n,\uu_n,v_n)$ satisfy the entropy  inequalities 
(\ref{EnergyEstimates}) and (\ref{BDEntropyEquality}) and converge to a weak solution of 
(\ref{TheSystemWithHydrostaticApproximation})-(\ref{BoundaryConditions}). 
\end{itemize}
\end{theo}
The proof of the main result is divided into three parts: the first part consists in writing System
(\ref{TheSystemWithHydrostaticApproximation}), using $(\xi,\uu, w=e^{-y} v)$ as unknowns instead of 
$(\rho,\uu, v)$. The obtained model is called \emph{model problem} (see Section
\ref{SectionModelProblem}). In the second part of the proof, we show the stability of weak solutions of the model
problem (see Section \ref{Step1}-\ref{Step5}). In the third and last part, by a simple criterion, the main result
is proved (see Section \ref{SectionEndOfTheProof}).
%%%%%%%%%%%%%%%%%%%%%%%%%%%%%%%%%%%%%%%%%%%%%%%%%%%%%%%%%%%%%%%%%%%%%%
% A simplified atmosphere model
%%%%%%%%%%%%%%%%%%%%%%%%%%%%%%%%%%%%%%%%%%%%%%%%%%%%%%%%%%%%%%%%%%%%%%
\section{Proof of the main result}\label{SubSectionChangeOfVariables}
The first part of the proof of Theorem \ref{MainRes} consists in writing the simplified model
(\ref{TheSystemWithHydrostaticApproximation}) in a more practical way, since the standard technique fails, as
pointed out in  Section \ref{SectionIntro}. 
\subsection{A model problem; an intermediate model}\label{SectionModelProblem}
We first begin, by noticing that the structure of the density  $\rho$ defined as a tensorial product (see
(\ref{DensityStratifie})) suggests the following change of variables:
\begin{equation}\label{ChangeOfVariables}
z = 1-e^{-y} 
\end{equation} where the vertical velocity,  in the new coordinates, is: 
\begin{equation}\label{DefinitionW}
w(t,x,z)=e^{-y}v(t,x,y). 
\end{equation}
Since the new vertical coordinate $z$ is defined as $\frac{d}{dy} z = e^{-y}$, 
multiplying by $e^{y}$ 
System (\ref{TheSystemWithHydrostaticApproximation}) and using the viscosity profile (\ref{Viscosity}) 
and the change of variables (\ref{ChangeOfVariables}) provides the following model:
\begin{equation}\label{SimplifiedAtmosphereModel}
\left\{\begin{array}{l}
\dsp \partial_t \xi + \mbox{div}_{x}\left(\xi\, \uu \right) + 
\partial_{z}\left(\xi \,w\right) =  0,\\
\dsp \partial_t \left(\xi\,\uu\right) +  
\mbox{div}_x\left(\xi\,\uu\otimes \uu\right) + 
\partial_z\left(\xi\,\uu\,w\right) +
\nabla_x \xi = \nu\mbox{div}_x\left(\xi D_x(\uu)\right) 
+\nu \partial_{zz}^2 (\xi \uu), \\
\dsp \partial_z \xi = 0  .
\end{array}\right. 
\end{equation}
which is the simplified CPEs (\ref{TheSystemWithHydrostaticApproximation})  with the unknowns
$$(\xi(t,x),\uu(t,x,y),w(t,x,y))\textrm{ instead of }(\rho(t,x,y),\uu(t,x,y),v(t,x,y))$$ that we call \emph{model
problem}.
In the new variables, the boundary conditions (\ref{BoundaryConditions}) and the initial conditions
(\ref{InitialCondition}) are written:
\begin{equation}\label{BoundaryConditionsForSimplifiedAtmosphereModel}
\begin{array}{l}
\textrm{ periodic conditions on }  \Omega_x,\\
w_{|z=0} = w_{|z=h} = 0,\\
{\partial_z \uu}_{|z=0} = {\partial_z \uu}_{|z=h}=0
\end{array}
\end{equation} and 
\begin{equation}\label{InitialConditionForSimplifiedAtmosphereModel}
\begin{array}{l}
 \uu(0,x,y) =\uu_0(x,z),\\
 \xi(0,x) =\xi_0(x)
 \end{array}
\end{equation} where  $\Omega = \mathbb{T}^2\times[0,h]$ with $h=1-e^{-1}$. 
%%%%%%%%%%%%%%%%%%%%%%%%%%%%%%%%%%%%%%%%%%%%%%%%%%%%%%%%%%%%%%%%%%%%%%
% Mathematical study of simplified atmosphere model
%%%%%%%%%%%%%%%%%%%%%%%%%%%%%%%%%%%%%%%%%%%%%%%%%%%%%%%%%%%%%%%%%%%%%%
\subsection{Mathematical study of the model problem}\label{SectionMathematicalStudy}
This section is devoted to the study of stability of weak solutions of  System 
(\ref{SimplifiedAtmosphereModel}) and equivalently for System  
(\ref{TheSystemWithHydrostaticApproximation}) as we will see in Section \ref{SectionEndOfTheProof}.
In what follows, we can say that:
\begin{definition}\label{DefinitionWeakSolution}
A weak solution of System (\ref{SimplifiedAtmosphereModel}) on $[0,T]\times\Omega$, with boundary 
(\ref{BoundaryConditionsForSimplifiedAtmosphereModel}) and 
initial conditions (\ref{InitialConditionForSimplifiedAtmosphereModel}), 
is a collection  of functions $(\xi,\uu,w)$, if 
$$
\begin{array}{ll}
\dsp \xi\in L^{\infty}(0,T;L^{3}({\Omega})),& \dsp \sqrt{\xi}\in L^{\infty}(0,T;H^{1}({\Omega})),\\
\dsp \sqrt{\xi}\uu\in L^2(0,T;L^2({\Omega})),& \dsp \sqrt{\xi}w\in L^\infty(0,T;(L^2({\Omega}))^2)\\
\dsp \sqrt{\xi}D_x(\uu)\in L^2(0,T;(L^2({\Omega}))^{2\times2}),& \sqrt{\xi}\partial_z w\in
L^2(0,T;L^2({\Omega})),\\
\dsp \nabla_x\sqrt{\xi}\in L^2(0,T;(L^2(\Omega))^2) &
\end{array}
$$
with  $\xi\geqslant 0$ and $(\xi,\sqrt{\xi}\uu,\sqrt{\xi}w)$ satisfies:
\begin{equation}\label{MassEquationDistribution}
\left\{
\begin{array}{l}
\dsp\partial_t \xi + \mbox{div}_x(\sqrt{\xi} \sqrt{\xi}\uu )+ \partial_z(\sqrt{\xi}\uu \sqrt{\xi}w)=0,\\
\xi(0,x) = \xi_0(x)
\end{array}
\right.
\end{equation} in the distribution sense, 
and the following equality holds for all smooth test function $\varphi$ 
with compact support such as $\varphi(T,x,z)=0$ 
and $\varphi_0=\varphi_{t=0}$:
\begin{equation}\label{def}
\begin{array}{l}
 -\int_0^T\int_{\Omega}\xi \uu\partial_t\varphi\,dx dz dt+
\int_0^T\int_{\Omega}\left(2\nu\xi D_x(\uu)-\xi \uu\otimes \uu\right):\nabla_x\varphi \,dx dz dt \\
-\int_0^T\int_{\Omega}\xi w\uu\partial_z\varphi \,dx dz dt
 -\nu\int_0^T\int_{\Omega}\xi \uu\partial_{zz}^2\varphi \,dx dz dt
+\int_0^T\int_{\Omega} r\xi\vert \uu\vert \uu\varphi\,dx dz dt \\
 -\int_0^T\int_{\Omega} \xi\mbox{div}_x(\varphi)\,dx dz dt=\int_{\Omega} \xi_{0}\uu_{0}\varphi_0 \,dx dz.
\end{array}
\end{equation}
\end{definition} 
We then have  the following result:
\begin{theo}\label{MainResult}
Let $(\xi_n,\uu_n,w_n)$ be a sequence of weak solutions of  
System (\ref{SimplifiedAtmosphereModel}), with boundary 
(\ref{BoundaryConditionsForSimplifiedAtmosphereModel}) and  
initial  conditions (\ref{InitialConditionForSimplifiedAtmosphereModel}),   satisfying entropy  inequalities 
(\ref{EnergyEstimates}) and (\ref{BDEntropyEquality}) such as 
\begin{equation}\label{ConvergenceInitialData}
\xi_n\geqslant 0,\quad \xi_0^n \to \xi_0 \textrm{ in } L^1({\Omega})
,\quad \xi_0^n \uu_0^n \to \xi_0 \uu_0 \textrm{ in } L^1({\Omega}). 
\end{equation}
Then, up to a subsequence,
\begin{itemize}
\item $\xi_n$ converges strongly in $\mathcal{C}^0(0,T; L^{3/2}({\Omega}))$,
\item $\sqrt{\xi_n} \uu_n$ converges strongly in $L^2(0,T; L^{3/2}({\Omega})^2)$,  
\item $\xi_n u_n$ converges strongly in $L^1(0,T; L^{1}({\Omega})^2)$ for all $T>0$,
\item $(\xi_n,\sqrt{\xi_n} \uu_n,\sqrt{\xi_n} w_n)$ converges to a weak solution 
of System (\ref{MassEquationDistribution}),
\item $(\xi_n,\uu_n,w_n)$ satisfy the entropy  inequalities 
(\ref{EnergyEstimates}) and (\ref{BDEntropyEquality}) and converge to a weak solution of 
(\ref{SimplifiedAtmosphereModel})-(\ref{BoundaryConditionsForSimplifiedAtmosphereModel}). 
\end{itemize}
\end{theo}
The proof of Theorem \ref{MainResult} is divided into three steps:
\begin{enumerate}
\item we first obtain suitable  \emph{a priori} bounds on $(\xi,\uu,w)$ (see Section \ref{SubSectionEstimates}). 
\item assuming the existence of sequences of weak solutions 
$(\xi_n,\uu_n,w_n)$, we show the compactness of sequences $(\xi_n,\uu_n,w_n)$ in apropriate space function (see
Section \ref{Step1}-\ref{Step4}).
\item using the obtained convergence, we show that we can pass to
the limit in all terms of System (\ref{SimplifiedAtmosphereModel}): this finishes the proof of Theorem
\ref{MainResult} (see Section \ref{Step5}).
\end{enumerate}
%%%%%%%%%%%%%%%%%%%%%%%%%%%%%%%%%%%%%%%%%%%%%%%%%%%%%%%%%%%%%%%%%%%%%%
% Energy and entropy estimates
%%%%%%%%%%%%%%%%%%%%%%%%%%%%%%%%%%%%%%%%%%%%%%%%%%%%%%%%%%%%%%%%%%%%%%
\subsubsection{Energy and entropy estimates}\label{SubSectionEstimates}
\noindent A part of \emph{a priori} bounds on  $(\xi,\uu,w)$ are obtained 
by the physical energy inequality which 
is obtained in a very classical way by multiplying the momentum equation by $\uu$, using
the mass equation and integrating by parts. We obtain the following  inequality:
\begin{equation}\label{EnergyEstimates}
\frac{d}{dt}\int_{\Omega}\big(\xi\frac{\uu^2}{2}+(\xi\ln\xi-\xi+1)\big)
+\int_{\Omega}\xi(\vert D_x(\uu)\vert^2+\vert\partial_z \uu\vert^2)
+r\int_{\Omega}\xi\vert \uu\vert^3\leqslant 0
\end{equation}
which provides  the uniform estimates:
\begin{eqnarray}
\sqrt{\xi} \uu\textrm{ is bounded in } L^\infty(0,T;(L^2({\Omega}))^2),\label{est2}\\
\xi^{\nicefrac{1}{3}}\uu\textrm{ is bounded in } L^3(0,T;(L^3({\Omega}))^2),\label{est3}\\
\sqrt{\xi}\partial_z\uu \textrm{ is bounded in } L^2(0,T;(L^2({\Omega}))^2),\label{est4}\\
\sqrt{\xi}D_x(\uu)\textrm{ is bounded in } L^2(0,T;(L^2({\Omega}))^{2\times 2}),\label{est5}\\
\xi\ln\xi-\xi+1\textrm{ is bounded in } L^\infty(0,T;L^1({\Omega})).\label{est6}
\end{eqnarray}
As pointed out by several authors (see e.g. \cite{BD03,MV07}), the crucial point in the proof of the stability in these
kind of models is to pass to the limit in the non linear term $\xi \uu\otimes \uu$ 
which requires the strong convergence of 
$\sqrt{\xi} \uu$. 
So we need additional information, which may be for instance provided by the mathematical BD-entropy \cite{BD02}: 

\noindent we first take the gradient of the mass equation, 
then we multiply by $2\nu $ and write the 
terms $\nabla_x \xi $ as $\xi \nabla_x \ln \xi$ to obtain:
\begin{equation}\label{Eq1}
\begin{array}{l}
\partial_t \left(2\nu \xi \nabla_x \ln\xi \right) + 
\mbox{div}_x \left(2\nu \xi \nabla_x \ln\xi \otimes \uu\right)
+ \partial_z \left(2\nu \xi \nabla_x \ln\xi w \right) \\ 
+ \mbox{div}_x \left(2\nu \xi \nabla_x^t \uu \right) 
+ \partial_z \left(2\nu \xi \nabla_x w \right) = 0 .
\end{array}
\end{equation} 
Next, we sum Equation (\ref{Eq1}) with the momentum 
equation of System (\ref{SimplifiedAtmosphereModel})
to get the equation:
\begin{equation}\label{Eq2}
\begin{array}{l}
\dsp \partial_t \left(\xi\,\psi\right) +  \mbox{div}_x\left(\psi \otimes \xi \uu\right) + 
\partial_z\left(\xi\,w\,\psi\right) + \partial_z\left(2\nu \xi\,\nabla_x w\right)  \\ 
+ \nabla_x \xi = 2\nu\mbox{div}_x \left(\xi A_x(\uu)\right)  - r \xi \abs{\uu}\uu 
+ \nu \xi \partial_z \left(\partial_z \uu\right)
\end{array}
\end{equation}
where $\psi = \uu + 2\nu \nabla_x \ln\xi$ and $\dsp 2 A_x(\uu) = \nabla_x \uu -\nabla_x^t \uu$
is the vorticity tensor.
The mathematical BD-entropy is then obtained by multiplying the previous equation by $\psi$ and integrating by
parts:
\begin{itemize}
 \item \begin{equation}\label{CalculBDEntropie1}
\int_{{\Omega}}\left(\partial_t \left(\xi\,\psi\right) 
+  \mbox{div}_x\left(\psi \otimes \xi \uu\right) + 
\partial_z\left(\xi\,w\,\psi\right)\right)\psi \,dx dz = \frac{d}{dt} 
\int_{{\Omega}} \xi \frac{\abs{\psi}^2}{2} \, dx dz.
 \end{equation} 
 \item Since $$\dsp \int_{{\Omega}} 2\nu\mbox{div}_x \left(\xi A_x(\uu)\right) 
\nabla_x \ln\xi \,dx dz = 0$$ and periodic boundary conditions are assumed on $\Omega_x$, we have:
\begin{equation}\label{CalculBDEntropie2}
\int_{{\Omega}}  2\nu\mbox{div}_x \left(\xi A_x(\uu)\right) \uu \,dx dz  =  
2\nu \int_{{\Omega}} \xi \abs{A_x(\uu)}^2\, dx dz.
 \end{equation} 
\item Derivating the mass equation with respect to $z$  provides the identity
$$\partial_z \mbox{div}_x (\xi \uu) = -\xi\partial_{zz}^2 w $$ and also recalling 
that $\xi$ is only $x$-dependent,
the integral becomes:
\begin{equation}\label{CalculBDEntropie3}
\begin{array}{lll}
\dsp \int_{{\Omega}}  \partial_z\left(2\nu \xi\,\nabla_x w\right) \psi  \,dx dz 
&=&\dsp  \int_{{\Omega}}\partial_z\left(2\nu \xi\,\nabla_x w\right) \uu  \,dx dz \\
&=&\dsp  \int_{{\Omega}}w \partial_z \mbox{div}_x (\xi \uu)  \,dx dz\\
&=&\dsp \int_{{\Omega}} \xi \abs{\partial_z w }^2  \,dx dz.
\end{array}
 \end{equation}
\item The other terms are easily computed. We have: 
\begin{equation}\label{CalculBDEntropie4}
\int_{{\Omega}}  r \xi \abs{\uu}\uu \psi \,dx dz = \int_{{\Omega}}  r \abs{\uu}^3 \,dx dz +
\int_{{\Omega}}  2\nu r \abs{\uu}\uu \nabla_x \xi \,dx dz,
\end{equation}
\begin{equation}\label{CalculBDEntropie5}
\int_{{\Omega}}  \nu \partial_z \left(\partial_z \uu\right) \psi \,dx dz = 
\int_{{\Omega}}  \nu \xi \abs{\partial_z \uu}^2 \,dx dz .
\end{equation}
and
\begin{equation}\label{CalculBDEntropie6}
\begin{array}{lll}
\dsp    \int_{{\Omega}} \nabla_x \xi \psi  \,dx dz 
&=&\dsp \int_{{\Omega}} \nabla_x\xi \uu \,dx dz + 
2\nu \int_{{\Omega}} \nabla_x\xi \nabla_x\ln\xi \,dx dz \\
&=&\dsp \frac{d}{dt}\int_{{\Omega}}(\xi\ln\xi -\xi +1) \,dx dz + 
8\nu \int_{{\Omega}} \abs{\nabla_x \sqrt{\xi}} \,dx dz.
\end{array}
 \end{equation}

\end{itemize}
Finally, gathering results (\ref{CalculBDEntropie1})--(\ref{CalculBDEntropie6})
provides the following mathematical entropy equality:
\begin{equation}\label{BDEntropyEquality}
\begin{array}{c}
\frac{d}{dt}\int_{\Omega}(\xi \frac{\abs{\psi}^2}{2}+ 
\xi\ln\xi-\xi+1) \,dx dz \\
+\int_{\Omega}(2\nu\xi\abs{\partial_z w}^2 
+ 2\nu \xi \abs{A_x(\uu)}^2 +\nu \xi\abs{\partial_{z} \uu}^2) \,dx dz\\
+\int_{\Omega} (r \xi \abs{\uu}^3+ 2\nu r \abs{\uu}\uu \nabla\xi
+ 8\nu \abs{\nabla_x\sqrt{\xi}}^2) \,dx dz= 0
\end{array}
\end{equation} and estimates:
\begin{eqnarray}
 \nabla\sqrt{\xi}\textrm{ is bounded in } L^\infty(0,T;(L^2({\Omega}))^3),\label{bd2}\\
% \xi^{\nicefrac{1}{3}}u\textrm{ is bounded in } L^3(0,T;(L^3({\Omega}))^2),\label{est3}\\
\sqrt{\xi}\partial_zw\textrm{ is bounded in } L^2(0,T;(L^2({\Omega}))),\label{bd4}\\
\sqrt{\xi}A_x(\uu)\textrm{ is bounded in } L^2(0,T;(L^2({\Omega}))^{2\times 2}).\label{bd5}
% (\xi\ln\xi-\xi+1)\big)\textrm{ is bounded in } L^\infty(0,T;L^1({\Omega})),\label{est6}
\end{eqnarray}
which finishes the first step of the proof of Theorem \ref{MainResult}.
\begin{rque}
Estimate (\ref{bd2}) is a straightforward consequence of estimates $\sqrt{\xi} \psi  \in  L^{\infty}(0,T,L^2(\Omega)^2)$ and 
$\sqrt{\xi} \uu  \in  L^{\infty}(0,T,L^2(\Omega)^2)$ since 
$$\sqrt{\xi}\psi=\sqrt{\xi}\uu + \frac{{\nabla_x \xi}}{\sqrt{\xi}}.$$ 
\end{rque}

The second step of the proof of Theorem \ref{MainResult}  will be divided into $4$ parts.
The first part consists to show the convergence of $\sqrt{\xi_n}$ (see Section \ref{Step1}). 
Then, we seek for bounds of 
$\sqrt{\xi_n} \uu_n$ and $\sqrt{\xi_n} w_n$ in an appropriate space (see Section \ref{Step2}) to be able to prove the
convergence 
of $\xi_n \uu_n$ (see Section \ref{Step3}). Thereafter, 
the convergence of $\sqrt{\xi_n} \uu_n$ is obtained (see Section \ref{Step4}).
%%%%%%%%%%%%%%%%%%%%%%%%%%%%%%%%%%%%%%%%%%%%%%%%%%%%%%%%%%%%%%%%%%%%%%
% Convergence of  $\sqrt{\xi_n}$
%%%%%%%%%%%%%%%%%%%%%%%%%%%%%%%%%%%%%%%%%%%%%%%%%%%%%%%%%%%%%%%%%%%%%%
\subsubsection{Convergence of  $\sqrt{\xi_n}$}\label{Step1}
The first part of the proof of Theorem \ref{MainResult} consists to show the following 
convergence.
\begin{lemme}\label{Lemma1}
For every $\xi_n$ satisfying the mass equation of System (\ref{SimplifiedAtmosphereModel}), we have: 
$$\sqrt{\xi_n}\textrm{ is bounded in }L^\infty(0,T,H^1({\Omega})),$$
$$\partial_t\sqrt{\xi_n}\textrm{ is bounded in }L^2(0,T,H^{-1}({\Omega})).$$
Then, up to a subsequence, the sequence $\xi_n$ converges almost everywhere
and strongly in $L^2(0,T ;L^2({\Omega}))$.
Moreover, $\xi_n$ converges to $\xi$ in $\mathcal{C}^0(0,T;L^{3/2}({\Omega}))$.
\end{lemme}
{\bf Proof of Lemma \ref{Lemma1}:}\\[0.5 \baselineskip]
The mass conservation equation gives
$$\vert\vert\sqrt{\xi_n}(t)\vert\vert^2_{L^2({\Omega})}
=\vert\vert\xi_{0}^n\vert\vert_{L^1({\Omega})}.$$
This equation and Estimate $(\ref{bd2})$ give the bound of $\sqrt{\xi_n}$ 
in $L^\infty(0,T,H^1({\Omega}))$.

\noindent Using again the mass conservation equation, we write
$$
\begin{array}{lll}
\dsp \partial_t(\sqrt{\xi_{n}})&=&\dsp -\frac{1}{2}\sqrt{\xi_{n}}\mbox{div}_x(\uu_{n})
-\uu_{n}.\nabla_x\sqrt{\xi_{n}}-\sqrt{\xi_n}\partial_zw_n\\
&=& \dsp\frac{1}{2}\sqrt{\xi_{n}}\mbox{div}_x(\uu_{n})
-\mbox{div}_x(\uu_{n}\sqrt{\xi_{n}})-\sqrt{\xi_n}\partial_zw_n.
\end{array}
$$
Then from Estimates (\ref{est5}), (\ref{bd5}), (\ref{bd2}), (\ref{bd4}),
$$\partial_t\sqrt{\xi_{n}} \textrm{ is bounded in } L^2(0,T,H^{-1}({\Omega})).$$
Next, Aubin's lemma provides compactness of $\sqrt{\xi_n}$ in  $\mathcal{C}^0(0,T,L^2({\Omega})$, that is:
\begin{equation}\label{cv1}
 \sqrt{\xi_{n}}\textrm{ converges strongly to }  
\sqrt{\xi}\textrm{ in }\mathcal{C}^0(0,T,L^2({\Omega})).
\end{equation}
We also have, by Sobolev embeddings, bounds of $\sqrt{\xi_n}$ in spaces 
$L^\infty(0,T,L^p({\Omega}))$ for all $p\in[1,6]$. \\Consequently, for $p=6$, we get bounds of  
$\xi_n$ in $L^\infty(0,T,L^3({\Omega}))$ and we deduce that
\begin{equation}\label{int1}
 \xi_n \uu_n=\sqrt{\xi_n}\sqrt{\xi_n}\uu_n\textrm{ is bounded in }L^\infty(0,T,L^{3/2}({\Omega})^2).
\end{equation}
It follows that $\partial_t\xi_n$ is bounded in $L^\infty(0,T,W^{-1,3/2}({\Omega}))$ since  
$$\partial_t \xi_n = -\mbox{div}(\xi_n \uu_n) - \xi_n\partial_z w_n$$ and we have Estimate (\ref{bd4}).

To conclude, writing $$\nabla_x\xi_n=2\sqrt{\xi_n}\nabla_x\sqrt{\xi_n}\in 
L^{\infty}(0,T; L^{3/2}({\Omega})^2),$$ we deduce bounds of $\xi_n$ in 
$L^{\infty}(0,T; W^{1,3/2}(\Omega))$.
%   By Sobolev embeddings, we have :
% $$W^{1,3/2}(\Omega)\subset\subset L^{3/2}(\Omega) \subset L^1(\Omega) \subset W^{-1,3/2}(\Omega)$$ 
% where $\subset\subset$ stands for compact embedding, and $\subset$ for continuous embedding.
Then Aubin's lemma provides  compactness of $\xi_n$ in the intermediate space $L^{3/2}({\Omega})$:
$$\textrm{ compactness  of } \xi_n  \textrm{ in } \mathcal{C}^0(0,T;L^{3/2}({\Omega})).$$
\findem
%%%%%%%%%%%%%%%%%%%%%%%%%%%%%%%%%%%%%%%%%%%%%%%%%%%%%%%%%%%%%%%%%%%%%%
% Bounds of  $\sqrt{\xi_n} w_n$ and $\sqrt{\xi_n} \uu_n$
%%%%%%%%%%%%%%%%%%%%%%%%%%%%%%%%%%%%%%%%%%%%%%%%%%%%%%%%%%%%%%%%%%%%%%
\subsubsection{Bounds of $\sqrt{\xi_n} \uu_n$ and $\sqrt{\xi_n} w_n$}\label{Step2}
To prove the convergence of the momentum equation, we have 
to control bounds of $\sqrt{\xi_n} \uu_n$ and $\sqrt{\xi_n} w_n$.
\begin{lemme}\label{Lemma2}
We have $$\sqrt{\xi_n} \uu_n \textrm{ bounded in } L^{\infty}(0,T;(L^2(\Omega))^2)$$
and 
$$\sqrt{\xi} w_n \textrm{ bounded in } L^2(0,T;L^2(\Omega)).$$
\end{lemme}
\noindent\textbf{Proof of Lemma \ref{Lemma2}:} 
We have already bounds of $\sqrt{\xi_n}$ (see Estimates (\ref{est2})). There is left to show bounds of 
$\sqrt{\xi_n} w_n$ $L^2(0,T;L^{2}(\Omega))$. 
\noindent As $\xi_n = \xi_n(t,x)$ and Estimates (\ref{bd4}) holds, by the Poincar\'e inequality, we have:  
$$\int_0^{h} \abs{\sqrt{\xi_n} w_n}^2 \,dz \leqslant c \int_0^{h}\abs{\partial_z (\sqrt{\xi_n} w_n)}^2\,dz.$$
Hence,  
$$\int_{\Omega}\xi_n \abs{ w_n}^2\,dx dz\leqslant c\int_{\Omega}\xi_n\abs{\partial_z w_n}^2\,dx dz$$
gives bounds of $\sqrt{\xi_n} w_n$  in $L^2(0,T;L^{2}({\Omega}))$.
\findem
%%%%%%%%%%%%%%%%%%%%%%%%%%%%%%%%%%%%%%%%%%%%%%%%%%%%%%%%%%%%%%%%%%%%%%
% Bounds of  $\sqrt{\xi_n} w_n$ and $\sqrt{\xi_n} \uu_n$
%%%%%%%%%%%%%%%%%%%%%%%%%%%%%%%%%%%%%%%%%%%%%%%%%%%%%%%%%%%%%%%%%%%%%%
\subsubsection{Convergence of $\xi_n \uu_n$}\label{Step3}
As bounds of $\sqrt{\xi_n} \uu_n$ and $\sqrt{\xi_n} w_n$ are provided by Lemma \ref{Lemma2}, we are able
to show the  convergence of the momentum equation.
\begin{lemme}\label{Lemma3}
Let $m_n=\xi_n \uu_n$ be a sequence satisfying the momentum 
equation (\ref{SimplifiedAtmosphereModel}). Then we have:
$$\xi_n\uu_n \rightarrow m \,\,\, \textrm{ in } L^2(0,T;(L^p(\Omega))^2) \textrm{ strong }, \forall\, 1\leqslant p<3/2$$ 
and
$$\xi_n\uu_n \rightarrow m \,\,\, \textrm{ a.e. } (t,x,y) \in (0,T)\times\Omega .$$ 
\end{lemme}
\noindent\textbf{Proof of Lemma \ref{Lemma3}:}\\
Writing  $\nabla_x(\xi_n \uu_n)$ as:
$$\nabla_x(\xi_n \uu_n)=\sqrt{\xi_n}\sqrt{\xi_n}\nabla_x \uu_n +2\sqrt{\xi_n} \uu_n\otimes\nabla\sqrt{\xi_n}$$ 
provides 
\begin{equation}\label{Step3Borne1}
\nabla_x(\xi_n \uu_n) \textrm{ bounded in } L^2(0,T;(L^{1}({\Omega}))^{2\times2}).
\end{equation}
Next, we have
\begin{equation}\label{Step3Borne2}
\partial_z(\xi_n \uu_n)=\sqrt{\xi_n}\sqrt{\xi_n}\partial_z(\uu_n)
\textrm{ is bounded } L^2(0,T;(L^{3/2}({\Omega}))^2).
\end{equation}
Then, from bounds (\ref{Step3Borne1}) and (\ref{Step3Borne2}), we deduce:
\begin{equation}\label{Step3Borne3}
\xi_n \uu_n\textrm{ is bounded } L^2(0,T;W^{1,1}({\Omega})^{2}).
\end{equation}
On the other hand, we have:
$$\begin{array}{lll}
\partial_t(\xi_n\,\uu_n)&=& -  
\mbox{div}_x\left(\xi_n\,\uu_n\otimes \uu_n\right) - 
\partial_z\left(\xi_n\,\uu_n\,w_n\right) -
\nabla_x \xi_n \\
& & + \nu\mbox{div}_x\left(\xi_n D_x(\uu_n)\right) 
+\nu \partial_z \left(\xi_n\partial_z \uu_n\right).
\end{array}
$$
As \begin{equation}\label{TermeNonLineaireDivergence}
   \xi_n \uu_n\otimes \uu_n = \sqrt{\xi}\uu_n \otimes  \sqrt{\xi}\uu_n,
\end{equation}
we deduce bounds of 
$$\xi_n \uu_n\otimes \uu_n\textrm{ in } L^\infty(0,T;(L^1({\Omega}))^{2\times2}).$$
% For all $ p>3$,  $W^{1,p}(\Omega)\subset L^{\infty}(\Omega)$, it follows that 
% $$L^1({\Omega})\subsetneq (L^{\infty}(\Omega))'\subset (W^{1,p}({\Omega}))' =W^{-1,p'}({\Omega})$$
% where $(X)'$ and $p'$ stand for the dual space of $X$, and the conjugate of $p$.
% Particularly, for $p=4$, we have $$\mbox{div}(\xi_n \uu_n\otimes \uu_n)\textrm{  bounded in  } 
Particularly, we have $$\mbox{div}(\xi_n \uu_n\otimes \uu_n)\textrm{  bounded in  } 
L^\infty(0,T;(W^{-2,4/3}({\Omega}))^{2}).$$
Similarly, as $\xi_n \uu_n  w_n = \sqrt{\xi}\uu_n  \sqrt{\xi}w_n \in (L^1(\Omega))^2$, we also have
$$\partial_z (\xi_n \uu_n  w_n)\textrm{  bounded in  } 
L^\infty(0,T;(W^{-2,4/3}({\Omega}))^{2}).$$
% \item  $\xi_n \uu_n  w_n = \sqrt{\xi}\uu_n  \sqrt{\xi}w_n \in (L^1(\Omega))^2,$
% \item  $\nabla_x \xi_n \in (L^1(\Omega))^2.$
% \end{itemize}
Moreover, as $$\sqrt{\xi_n}\sqrt{\xi_n}\partial_z \uu_n, \in L^2(0,T;(L^{3/2}({\Omega}))^2) \textrm{ and}$$ 
$$\sqrt{\xi_n}\sqrt{\xi_n}D_x(\uu_n)\in L^2(0,T;(L^{3/2}({\Omega}))^{2\times2}), $$ 
we get bounds of 
$$\partial_z(\sqrt{\xi_n}\sqrt{\xi_n}\partial_z \uu_n),  
\,\mbox{div}_x(\sqrt{\xi_n}\sqrt{\xi_n}D_x(\uu_n))
\in L^2(0,T;(W^{-1,3/2}({\Omega}))^2).$$
We also have  bounds of $\nabla_x\xi_n \in L^\infty(0,T,(W^{-1,3/2}({\Omega}))^2).$\\
Using 
% $L^1({\Omega})\subset W^{-1,4/3}({\Omega})$ and 
$W^{-1,3/2}({\Omega})\subset W^{-1,4/3}({\Omega})$, we obtain 
\begin{equation}\label{Step3Borne4}
 \partial_t(\xi_n \uu_n)\textrm{ bounded in }L^2(0,T; W^{-2,4/3}({\Omega})^2).
\end{equation}
Using bounds (\ref{Step3Borne3}), (\ref{Step3Borne4}) with  Aubin's lemma provides compactness of
\begin{equation}\label{int5}
\xi_n \uu_n\in L^2(0,T; (L^{p}({\Omega}))^2),\,\forall p\in\,[1,3/2[.
\end{equation}
\findem
%%%%%%%%%%%%%%%%%%%%%%%%%%%%%%%%%%%%%%%%%%%%%%%%%%%%%%%%%%%%%%%%%%%%%%
% Bounds of  $\sqrt{\xi_n} w_n$ and $\sqrt{\xi_n} \uu_n$
%%%%%%%%%%%%%%%%%%%%%%%%%%%%%%%%%%%%%%%%%%%%%%%%%%%%%%%%%%%%%%%%%%%%%%
\subsubsection{Convergence of $\sqrt{\xi_n} \uu_n$ and ${\xi_n} w_n$}\label{Step4}
Let us note that, up to Section \ref{Step3}, we can always 
define $\uu = m/\xi$ on the set $\{\xi>0\}$, but we do not know, \emph{a priori}, if $m$ equals  zero  
on the vacuum set. To this end, we need to prove the following lemma:
\begin{lemme}\label{Lemma4}
\begin{enumerate}
\item[]
 \item The sequence $\sqrt{\xi_{n}} \uu_{n}$ satisfies
\begin{itemize}
 \item $\sqrt{\xi_{n}} \uu_{n}$ converges strongly 
in $L^2(0,T;L^2({\Omega}))$ to $\frac{m}{\sqrt{\xi}}.$
 \item We have $m=0$ almost everywhere on the set $\{\xi=0\}$ 
and there exists a function  $\uu$ such that  $m=\xi \uu$ and
\begin{equation}\label{cv3}
 \xi_{n} \uu_{n}\to \xi \uu\textrm{ strongly in }
L^2(0,T;L^p({\Omega})^2)\textrm{ for all }p\in [1,3/2[,
\end{equation}
\begin{equation}\label{cv4}
 \sqrt{\xi_{n}} \uu_{n}\to \sqrt{\xi} \uu\textrm{  strongly in  }L^2(0,T;L^2({\Omega})^2).
\end{equation}
\end{itemize}
\item The sequence $\sqrt{\xi_{n}} w_{n}$  
converges weakly in $L^2(0,T;L^2({\Omega}))$ to $\sqrt{\xi} w$.
\end{enumerate}
\end{lemme}
\textbf{Proof of Lemma \ref{Lemma4}:}\\
We refer to \cite{NS10} for details of the first part of the proof.
The second part of the theorem is done by weak compactness. As $\sqrt{\xi_n} w_n$ is bounded
in $L^2(0,T;L^2({\Omega}))$,  there exists, up to a subsequence, 
$\sqrt{\xi_n} w_n$ which converges weakly some limit $l$ in  $L^2(0,T;L^2({\Omega}))$. Next,  we define  $w$ 
$$
w = \left\{
\begin{array}{lll}
\dsp \frac{l}{\sqrt{\xi}} & \textrm{ if }& \xi>0, \\
0                         \textrm{ a.e. }& \textrm{ if }& \xi=0\,
\end{array}\right.
$$
where the limit $l$ is written: $\dsp l = \sqrt{\xi}\frac{l}{\sqrt{\xi}}  = \sqrt{ \xi} w$.
\findem
This finishes the second step of the proof of Theorem \ref{MainResult}.

In the third and last step (see Section \ref{Step5}), 
using the convergence of Sections \ref{Step1}--\ref{Step4}, we show that we can pass to the limit for all
terms of System (\ref{SimplifiedAtmosphereModel}).
%%%%%%%%%%%%%%%%%%%%%%%%%%%%%%%%%%%%%%%%%%%%%%%%%%%%%%%%%%%%%%%%%%%%%%
% Convergence step 
%%%%%%%%%%%%%%%%%%%%%%%%%%%%%%%%%%%%%%%%%%%%%%%%%%%%%%%%%%%%%%%%%%%%%%
\subsubsection{Convergence step}\label{Step5}
We are now ready to prove that we can pass to the limit in all terms of System (\ref{SimplifiedAtmosphereModel}) 
in the sense of Theorem \ref{MainResult}.
To this end, let ($\xi_n$, $\uu_n$, $w_n$) be a weak solution of System 
(\ref{SimplifiedAtmosphereModel}) satisfying Lemma \ref{Lemma1} to \ref{Lemma4} and 
let $\phi \in \mathcal{C}_c^{\infty}([0,T]\times\Omega)$ be a smooth function with compact support 
such as $\phi(T,x,z)=0$. Then, writing each term of the weak formulation of 
System (\ref{SimplifiedAtmosphereModel}), we have:
\begin{itemize}
\item
\begin{equation}
\begin{array}{lll}
\int_0^T \int_\Omega \partial_t (\xi_n \uu_n)\phi \,dx dz dt &=&-\int_0^T\int_\Omega \xi_n \uu_n
\partial_t \phi \,dx dz dt \\
& &  -\int_\Omega \xi_0^n \uu_0^n \phi(0,x,z)\,dx dz.
\end{array}
\end{equation}

Using convergences (\ref{ConvergenceInitialData}) and Lemma \ref{Lemma3}, we get
$$-\int_0^T\int_\Omega \xi_n \uu_n
\partial_t \phi \,dx dz dt-\int_\Omega \xi_0^n \uu_0^n \phi(0,x,z)\,dx dz \to $$ 
$$-\int_0^T\int_\Omega\xi \uu \partial_t \phi \,dx dz dt-\int_\Omega \xi_0 \uu_0 \phi(0,x,y)\,dx dz.$$
\item 
$$\int_0^T\int_\Omega\mbox{div}_x(\xi_n \uu_n\otimes \uu_n)\cdot\phi \,dx dz dt
=-\int_0^T\int_\Omega \xi_n \uu_n\otimes \uu_n:\nabla_x\phi\,dx dz dt.$$
From Equality (\ref{TermeNonLineaireDivergence}) and Lemma \ref{Lemma4}, we have: 
$$-\int_0^T\int_\Omega \xi_n \uu_n\otimes \uu_n:\nabla_x\phi\,dx dz dt 
\to -\int_0^T\int_\Omega \xi \uu\otimes \uu:\nabla_x\phi\,dx dz dt.$$
\item 
$$ \int_0^T\int_\Omega \partial_z(\xi_n \uu_n w_n)\cdot \phi \, dx dz dt=
-\int_0^T\int_\Omega \xi_n \uu_n w_n \cdot \partial_z \phi \, dx dz dt.
$$
As $\xi_n \uu_n w_n = \sqrt{\xi_n}\uu_n \sqrt{\xi_n}w_n $, by Lemma \ref{Lemma4}, we get:
$$ -\int_0^T\int_\Omega \xi_n \uu_n w_n \cdot \partial_z \phi  \, dx dz dt \to
-\int_0^T\int_\Omega \xi \uu w \cdot \partial_z \phi \, dx dz dt.
$$
\item
$$ \int_0^T\int_\Omega \nabla_x \xi_n \cdot \phi\, dx dz dt=
-\int_0^T\int_\Omega \xi_n \mbox{div}_x(\phi)\, dx dz dt.
$$
Then, Lemma \ref{Lemma1} provides:
$$
-\int_0^T\int_\Omega \xi_n \mbox{div}_x(\phi)\, dx dz dt
\to
-\int_0^T\int_\Omega \xi \mbox{div}_x(\phi)\, dx dz dt
$$
\item
$$\int_0^T\int_\Omega\mbox{div}_x(\xi_n D_x(\uu_n))\cdot\phi\,dx dz dt
=-\int_0^T\int_\Omega \xi_n D_x(\uu_n):\nabla\phi\,dx dz dt.$$
Since $D_x(\uu_n)=\frac{1}{2}(\nabla_x \uu_n+\nabla_x^t \uu_n)$, expanding the term in the 
last integral gives:
$$
\begin{array}{l}
-\dsp \int_0^T\int_\Omega \xi_n D_x(\uu_n):\nabla_x\phi\, dx dz dt \\ = 
\dsp \frac{1}{2}\int_0^T\int_\Omega( \xi_n\uu_n\cdot\Delta_x\phi
+\nabla_x\phi \nabla_x(\sqrt{\xi_n})\cdot \sqrt{\xi_n}\uu_n\, dx dz dt \\
\dsp + \frac{1}{2}\int_0^T\int_\Omega( \xi_n\uu_n\cdot\mbox{div}_x(\nabla_x^t \phi)
+\nabla_x^t\sqrt{\xi_n}\cdot\nabla_x\phi\cdot\sqrt{\xi_n}\uu_n)\, dx dz dt.
\end{array}
$$
From Estimates (\ref{bd2}), the sequence $\nabla_x \sqrt{\xi_n}$ weakly converges,
and using Lemma \ref{Lemma1}, Lemma \ref{Lemma3} and \ref{Lemma4}, we obtain:
$$\dsp \frac{1}{2}\int_0^T\int_\Omega( \xi_n\uu_n\cdot\Delta_x\phi
+\nabla_x\phi \nabla_x(\sqrt{\xi_n})\cdot \sqrt{\xi_n}\uu_n\, dx dz dt $$
$$\dsp + \frac{1}{2}\int_0^T\int_\Omega( \xi_n\uu_n\cdot\mbox{div}_x(\nabla_x^t \phi)
+\nabla_x^t\sqrt{\xi_n}\cdot\nabla_x\phi\cdot\sqrt{\xi_n}\uu_n)\, dx dz dt\to$$
$$\dsp \frac{1}{2}\int_0^T\int_\Omega( \xi \uu \cdot\Delta_x\phi
+\nabla_x\phi \nabla_x(\sqrt{\xi})\cdot \sqrt{\xi}\uu\, dx dz dt $$
$$\dsp + \frac{1}{2}\int_0^T\int_\Omega( \xi \uu \cdot\mbox{div}_x(\nabla_x^t \phi)
+\nabla_x^t\sqrt{\xi}\cdot\nabla_x\phi\cdot\sqrt{\xi}\uu)\, dx dz dt.$$
Hence
$$-\int_0^T\int_\Omega \xi_n D_x(\uu_n):\nabla_x\phi\,dx dz dt\to 
-\int_0^T\int_\Omega \xi D_x(\uu):\nabla_x\phi\,dx dz dt.$$
\item $$\int_0^T\int_\Omega \partial_{zz}^2(\xi_n \uu_n) \cdot\phi\,dx dz dt\to 
\int_0^T\int_\Omega \xi_n \uu_n \cdot\partial_{zz}^2(\phi)\,dx dz dt.$$
Using Lemma \ref{Lemma3} provides the following convergence:
$$\int_0^T\int_\Omega \xi_n \uu_n \cdot\partial_{zz}^2(\phi)\,dx dz dt \to 
\int_0^T\int_\Omega \xi \uu \cdot\partial_{zz}^2(\phi)\,dx dz dt
$$
\item $$\int_0^T\int_{\Omega} r\xi_n \abs{\uu_n} \uu_n \cdot\phi\,dx dz dt \to 
\int_0^T\int_{\Omega} r\xi \abs{\uu} \uu \cdot\phi\,dx dz dt 
$$ with Lemma \ref{Lemma4}, which finishes the proof of Theorem  \ref{MainResult}.
\end{itemize}
\findem
%%%%%%%%%%%%%%%%%%%%%%%%%%%%%%%%%%%%%%%%%%%%%%%%%%%%%%%%%%%%%%%%%%%%%%
%  Proof of Corollary \ref{Corollaire}
%%%%%%%%%%%%%%%%%%%%%%%%%%%%%%%%%%%%%%%%%%%%%%%%%%%%%%%%%%%%%%%%%%%%%%
\subsubsection{End of the proof of Theorem \ref{MainRes}}\label{SectionEndOfTheProof}
In order to conlude, let ($\xi_n$, $\uu_n$, $w_n$) be a weak solution of System 
(\ref{SimplifiedAtmosphereModel}), then all estimates \ref{Step1}-\ref{Step5}
hold if we replace $\xi_n$ by $\rho_n$ and $w_n$ by $v_n$ (see \cite{EN10_3}), since $\rho(t,x,y)=\xi(t,x)
e^{-y}$ and $w(t,x,z) = v(t,xy) e^{-y}$ where $\frac{d}{dy}z = e^{-y}$. This proves Theorem \ref{MainRes}.
\findem
%%%%%%%%%%%%%%%%%%%%%%%%%%%%%%%%%%%%%%%%%%%%%%%%%%%%%%%%%%%%%%%%%%%%%%
%  Biblio
%%%%%%%%%%%%%%%%%%%%%%%%%%%%%%%%%%%%%%%%%%%%%%%%%%%%%%%%%%%%%%%%%%%%%%


\begin{thebibliography}{10}

\bibitem{HandBookBresch09}
D.~Bresch.
\newblock Shallow-water equations and related topics.
\newblock In {\em Handbook of differential equations: evolutionary equations.
  {V}ol. {V}}, Handb. Differ. Equ., pages 1--104. Elsevier/North-Holland,
  Amsterdam, 2009.

\bibitem{BD02}
D.~Bresch and B.~Desjardins.
\newblock Sur un mod\`ele de {S}aint-{V}enant visqueux et sa limite
  quasi-g\'eostrophique.
\newblock {\em C. R. Math. Acad. Sci. Paris}, 335(12):1079--1084, 2002.

\bibitem{BD03}
D.~Bresch and B.~Desjardins.
\newblock Existence of global weak solutions for a 2{D} viscous shallow water
  equations and convergence to the quasi-geostrophic model.
\newblock {\em Comm. Math. Phys.}, 238(1-2):211--223, 2003.

\bibitem{BD04}
D.~Bresch and B.~Desjardins.
\newblock Some diffusive capillary models of {K}orteweg type.
\newblock {\em C. R. Acad. Sciences}, 332(11):881--886, 2004.

\bibitem{BD06}
D.~Bresch and B.~Desjardins.
\newblock Sur la th\'eorie globale des \'equations de {N}avier-{S}tokes
  compressible.
\newblock {\em Journ\'ees \'equations aux d\'eriv\'ees partielles}, Exp. No. 3,
  26 p, 2006.

\bibitem{BDGV07}
D.~Bresch, B.~Desjardins, and D.~G\'erard-Varet.
\newblock On compressible {N}avier-{S}tokes equations with density dependent
  viscosities in bounded domains.
\newblock {\em J. Math. Pures Appl. (9)}, 87(2):227--235, 2007.

\bibitem{BDGG07}
D.~Bresch, B.~Desjardins, J.-M. Ghidaglia, and E.~Grenier.
\newblock Mathematical properties of the basic two fluid model.
\newblock {\em to appear in Arch. Rat. Mech. Anal.}, 2007.

\bibitem{B81}
G.~Buntebarth.
\newblock Zur entwicklung des begriffes geophysik.
\newblock {\em Abhandlungen der Braunschwiegischen Wissenschaftlichen
  Gesellschaft}, 32:95--109, 1981.

\bibitem{EN10_3}
M.~Ersoy and T.~Ngom.
\newblock Existence of a global weak solution to one model of compressible
  primitive equations.
\newblock {\em Submitted}, 2010.

\bibitem{GK05}
B.~V. Gatapov and A.~V. Kazhikhov.
\newblock Existence of a global solution to one model problem of atmosphere
  dynamics.
\newblock {\em Siberian Mathematical Journal}, 46(5):805--812, 2005.

\bibitem{K36}
N.~E. Kochin.
\newblock On simplification of the equations of hydromechanics in the case of
  the general circulation of the atmosphere.
\newblock {\em Trudy Glavn. Geofiz. Observator.}, 4:21--45, 1936.

\bibitem{LTW92}
J.L. Lions, R.~Temam, and S.~Wang.
\newblock New formulations for the primitive equations for the atmosphere and
  applications.
\newblock {\em Nonlinearity}, 5:237--288, 1992.

\bibitem{MV07}
A.~Mellet and A.~Vasseur.
\newblock On the barotropic compressible navier-stokes equations.
\newblock {\em Comm. Partial Differential Equations}, 32(1-3):431--452, 2007.

\bibitem{NS10}
T.~Ngom and M.~Sy.
\newblock Derivation and stability study of a rigid lid bilayer model.
\newblock {\em Submitted}, 2010.

\bibitem{Pedlowski87}
J.~Pedlowski.
\newblock {\em Geophysical Fluid dynamics}.
\newblock 2nd Edition, Springer-Verlag, New-York, 1987.

\bibitem{Temam77}
R.~Temam.
\newblock Navier-{S}tokes equations.
\newblock pages xiv+408, 2001.
\newblock Theory and numerical analysis, Reprint of the 1984 edition.

\bibitem{TZ04}
R.~Temam and M.~Ziane.
\newblock {\em Some mathematical problems in geophysical fluid dynamics}.
\newblock North-Holland, Amsterdam, 2004.

\end{thebibliography}
\end{document}